\documentclass[11pt]{amsart}

\usepackage[all]{xy}
\usepackage{graphics}

\catcode`\Ž=\active \def Ž{\'e}
\catcode`\=\active \def {\`e}
\catcode`\ˆ=\active \def ˆ{\`a} 
\catcode`\=\active \def {\`u}

\catcode`\=\active \def {\^e}        
\catcode`\"=\active \def "{\^{\i}}    
\catcode`\â=\active \def â{\^a}

\catcode`\"=\active \def "{\^{\i}}
\catcode`\ô=\active \def ô{\^o}
\catcode`\ž=\active \def ž{\^u}
\catcode`\é=\active \def é{\`E} 
\catcode`\æ=\active \def æ{\^E}

\newtheorem{lemma}[subsection]{Lemma}

\catcode`\Ž=\active \def Ž{\'e}
\catcode`\=\active \def {\`e}
\catcode`\ˆ=\active \def ˆ{\`a} 
\catcode`\=\active \def {\`u}

\catcode`\=\active \def {\^e}        
\catcode`\"=\active \def "{\^{\i}}    
\catcode`\â=\active \def â{\^a}

\catcode`\"=\active \def "{\^{\i}}
\catcode`\ô=\active \def ô{\^o}
\catcode`\ž=\active \def ž{\^u}
\catcode`\é=\active \def é{\`E} 
\catcode`\æ=\active \def æ{\^E}

\def\KK{{\mathcal K}}

\newenvironment{proo}{\begin{trivlist} \item{\emph{Proof.}}}
  {\hfill $\square$ \end{trivlist}}

\def\KK{{\mathbb{K}}}

\def\gg{\mathfrak{g}}

\def\arbreA{\vcenter{\xymatrix@R=3pt@C=3pt{
&& \\
&*{}\ar@{-}[ur] \ar@{-}[ul] \ar@{-}[d]     &\\
&&
}}}

\def\arbreAgrand{\vcenter{\xymatrix@R=30pt@C=30pt{
&& \\
&*{}\ar@{-}[ur] \ar@{-}[ul] \ar@{-}[d]     &\\
&&
}}}

\def\arbreBA{\vcenter{\xymatrix@R=2pt@C=2pt{
&&&&\\
&&&*{}\ar@{-}[ul] & \\
&&*{}\ar@{-}[uurr] \ar@{-}[uull] \ar@{-}[d]     &&\\
&&&&
}}}

\def\arbreBAlabelled{\vcenter{\xymatrix@R=6pt@C=6pt{
1&&0&&1\\
&&&*{}\ar@{-}[ul] & \\
&&*{}\ar@{-}[uurr] \ar@{-}[uull] \ar@{-}[d]     &2&\\
&&0&&
}}}

\def\arbreAB{\vcenter{\xymatrix@R=2pt@C=2pt{
&&&&\\
&*{}\ar@{-}[ur] &&& \\
&&*{}\ar@{-}[uurr] \ar@{-}[uull] \ar@{-}[d]     &&\\
&&&&
}}}

\def\arbreABlabelled{\vcenter{\xymatrix@R=6pt@C=6pt{
1&&0&&1\\
&*{}\ar@{-}[ur] &&& \\
&2&*{}\ar@{-}[uurr] \ar@{-}[uull] \ar@{-}[d]     &&\\
&&0&&
}}}

\def\arbreBB{\vcenter{\xymatrix@R=2pt@C=2pt{
&&*{}&&\\
&&&& \\
&&*{}\ar@{-}[uurr] \ar@{-}[uull] \ar@{-}[d] \ar@{-}[uu]     &&\\
&&&&
}}}

\def\arbreBBbis{\vcenter{\xymatrix@R=4pt@C=4pt@M=0pt{
&&*{}&&\\
&&&& \\
&&*{}\ar@{-}[uurr] \ar@{-}[uull] \ar@{-}[d] \ar@{-}[uu]     &&\\
&&&& }}}

\def\arbreABC{\vcenter{\xymatrix@R=1pt@C=0pt{
&&&&&&\\
&*{}\ar@{-}[ur] &&&&& \\
&&*{}\ar@{-}[uurr] &&&&\\
&&&*{}\ar@{-}[uuurrr] \ar@{-}[uuulll] \ar@{-}[d] &&&\\
&&&&&&
}}}

\def\arbreABCplus{\vcenter{\xymatrix@R=1pt@C=0pt{
&&&&&&\\
*{}\ar@{-}[u]&&*{}\ar@{-}[u]&&*{}\ar@{-}[u]&&*{}\ar@{-}[u]\\
&*{}\ar@{-}[ur] &&&&& \\
&&*{}\ar@{-}[uurr] &&&&\\
&&&*{}\ar@{-}[uuurrr] \ar@{-}[uuulll] \ar@{-}[d] &&&\\
&&&&&&
}}}

\def\arbreABClabelled{\vcenter{\xymatrix@R=5pt@C=5pt{
1&&2&&1&&1\\
&*{}\ar@{-}[ur] &&&&& \\
&0&*{}\ar@{-}[uurr] &&&&\\
&&2&*{}\ar@{-}[uuurrr] \ar@{-}[uuulll] \ar@{-}[d] &&&\\
&&&0&&&
}}}

\def\arbreBAC{\vcenter{\xymatrix@R=1pt@C=0pt{
&&&&&&\\
&&&*{}\ar@{-}[ul] &&& \\
&&*{}\ar@{-}[uurr] &&&&\\
&&&*{}\ar@{-}[uuurrr] \ar@{-}[uuulll] \ar@{-}[d] &&&\\
&&&&&&
}}}

\def\arbreACA{\vcenter{\xymatrix@R=1pt@C=0pt{
&&&&&&\\
&*{}\ar@{-}[ur] &&&&*{}\ar@{-}[ul] & \\
&&&&&&\\
&&&*{}\ar@{-}[uuurrr] \ar@{-}[uuulll] \ar@{-}[d] &&&\\
&&&&&&
}}}

\def\arbreCAB{\vcenter{\xymatrix@R=1pt@C=0pt{
&&&&&&\\
&&&*{}\ar@{-}[ur] &&& \\
&&&&*{}\ar@{-}[uull] &&\\
&&&*{}\ar@{-}[uuurrr] \ar@{-}[uuulll] \ar@{-}[d] &&&\\
&&&&&&
}}}

\def\arbreCABplus{\vcenter{\xymatrix@R=1pt@C=0pt{
&&&&&&\\
*{}\ar@{-}[u]&&*{}\ar@{-}[u]&&*{}\ar@{-}[u]&&*{}\ar@{-}[u]\\
&&&*{}\ar@{-}[ur] &&& \\
&&&&*{}\ar@{-}[uull] &&\\
&&&*{}\ar@{-}[uuurrr] \ar@{-}[uuulll] \ar@{-}[d] &&&\\
&&&&&&
}}}

\def\arbreCABlabelled{\vcenter{\xymatrix@R=5pt@C=5pt{
1&&2&&1&&1\\
&&&*{}\ar@{-}[ur] &&& \\
&&&0&*{}\ar@{-}[uull] &&\\
&&&*{}\ar@{-}[uuurrr] \ar@{-}[uuulll] \ar@{-}[d] &2&&\\
&&&0&&&
}}}

\def\arbreCBA{\vcenter{\xymatrix@R=1pt@C=0pt{
&&&&&&\\
&&&&&*{}\ar@{-}[ul] & \\
&&&&*{}\ar@{-}[uull] &&\\
&&&*{}\ar@{-}[uuurrr] \ar@{-}[uuulll] \ar@{-}[d] &&&\\
&&&&&&
}}}

\def\arbreAAC{\vcenter{\xymatrix@R=1pt@C=0pt{
&&&&&&\\
&&&&&& \\
&&*{}\ar@{-}[uurr]\ar@{-}[uu]  &&&&\\
&&&*{}\ar@{-}[uuurrr] \ar@{-}[uuulll] \ar@{-}[d] &&&\\
&&&&&&
}}}

\def\arbreCAC{\vcenter{\xymatrix@R=1pt@C=0pt{
&&&&&&\\
&&&*{}\ar@{-}[ul]\ar@{-}[ur] &&& \\
&&&&&&\\
&&&*{}\ar@{-}[uuurrr] \ar@{-}[uuulll]\ar@{-}[uu]  \ar@{-}[d] &&&\\
&&&&&&
}}}

\def\arbreACC{\vcenter{\xymatrix@R=1pt@C=0pt{
&&&&&&\\
&*{}\ar@{-}[ur] &&&&& \\
&&&&&&\\
&&&*{}\ar@{-}[uuurrr] \ar@{-}[uuulll]\ar@{-}[uuu]  \ar@{-}[d] &&&\\
&&&&&&
}}}

\def\arbreCBB{\vcenter{\xymatrix@R=1pt@C=0pt{
&&&&&&\\
&&&&&& \\
&&&&*{}\ar@{-}[uull]\ar@{-}[uu]  &&\\
&&&*{}\ar@{-}[uuurrr] \ar@{-}[uuulll] \ar@{-}[d] &&&\\
&&&&&&
}}}

\def\arbreCCA{\vcenter{\xymatrix@R=1pt@C=0pt{
&&&&&&\\
&&&&&*{}\ar@{-}[ul] & \\
&&&&&&\\
&&&*{}\ar@{-}[uuurrr] \ar@{-}[uuulll]\ar@{-}[uuu]  \ar@{-}[d] &&&\\
&&&&&&
}}}

\def\arbreBBC{\vcenter{\xymatrix@R=1pt@C=0pt{
&&&&&&\\
&&&&&& \\
&&*{}\ar@{-}[uu] \ar@{-}[uurr] &&&&\\
&&&*{}\ar@{-}[uuurrr] \ar@{-}[uuulll] \ar@{-}[d] &&&\\
&&&&&&
}}}

\def\arbreCCC{\vcenter{\xymatrix@R=1pt@C=0pt{
&&&&&&\\
&&&&&& \\
&&&&&&\\
&&&*{}\ar@{-}[uuurrr]  \ar@{-}[uuulll]\ar@{-}[uuul] \ar@{-}[uuur]  \ar@{-}[d] &&&\\
&&&&&&
}}}

\def\arbreBAEAB{\vcenter{\xymatrix@R=4pt@C=3pt{
&&&& &&&& &&\\
*{}\ar@{-}[u]&&*{}\ar@{-}[u]& &*{}\ar@{-}[u]&&*{}\ar@{-}[u]&&*{}\ar@{-}[u]&&*{}\ar@{-}[u]\\
&& & *{}\ar@{-}[ul]& &&& *{}\ar@{-}[ur] & && \\
&& *{}\ar@{-}[uurr]&& &&&& *{}\ar@{-}[uull]&&\\
&&&& &&&& &&\\
&&&& &&&& &&\\
&&&& &*{}\ar@{-}[uuuuulllll]\ar@{-}[uuuuurrrrr] \ar@{-}[d]&&& &&\\
&&&& &&&& &&
}}}

\def\arbreAEACA{\vcenter{\xymatrix@R=4pt@C=3pt{
&&&& &&&& &&\\
*{}\ar@{-}[u]&&*{}\ar@{-}[u]& &*{}\ar@{-}[u]&&*{}\ar@{-}[u]&&*{}\ar@{-}[u]&&*{}\ar@{-}[u]\\
& *{}\ar@{-}[ur]& && & *{}\ar@{-}[ur] &&&& *{}\ar@{-}[ul]  & \\
&&&& &&& &&&\\
&&&& &&&*{}\ar@{-}[uuulll]& &&\\
&&&& &&&& &&\\
&&&& &*{}\ar@{-}[uuuuulllll]\ar@{-}[uuuuurrrrr] \ar@{-}[d]&&& &&\\
&&&& &&&& &&
}}}

\def\graphUN
{\xymatrix@R=12pt@C=12pt{
&&&&*{}\ar@{-}[r] &*{}\ar@{-}[dddrrr]&&& & \\
&&&*{}\ar@{-}[dr]& &&&& & \\
&&*{}\ar@{-}[ddrr]&&*{}\ar@{-}[r] &*{}\ar@{-}[dr]&&& & \\
*{}\ar@{-}[r]&*{}\ar@{-}[uuurrr]\ar@{-}[dddrrr]&&& &&*{}\ar@{-}[dr]&&*{}\ar@{-}[r] & \\
&&&&*{}\ar@{-}[r] &*{}\ar@{-}[ur]&&*{}& & \\
&&&& &&&& & \\
&&&&*{}\ar@{-}[r] &*{}\ar@{-}[uuurrr]&&& & \\
}}

\def\graphUNsol
{\xymatrix@R=12pt@C=12pt{
&&&&  &&&& & \\
&&&&*{}\ar@{-}[r]^-{1} &*{}\ar@{-}[dddrrr]&&& & \\
&&&*{}\ar@{-}[dr]\ar@{-}[ur]& &&&& & \\
&&*{}\ar@{-}[ddrr]\ar@{-}[ur]^-{0}&&*{}\ar@{-}[r]^-{2} &*{}\ar@{-}[dr]&&& & \\
*{}\ar@{-}[r]^-{0}&*{}\ar@{-}[ur]^-{2}\ar@{-}[dddrrr]&&& &&*{}\ar@{-}[dr]_-{0}&&*{}\ar@{-}[r]^-{0} & \\
&&&&*{}\ar@{-}[r]^-{1} &*{}\ar@{-}[ur]&&*{}\ar@{-}[ur]_-{2}& & \\
&&&& &&&& & \\
&&&&*{}\ar@{-}[r]^-{1} &*{}\ar@{-}[uurr]&&& & \\
}}

\begin{document}

\title{The YY game}
\author{Jean-Louis Loday}
\address{Institut de Recherche Math\'ematique Avanc\'ee\\
    CNRS et Universit\'e de Strasbourg\\
    7 rue R. Descartes,
    67084 Strasbourg Cedex, France}
    
\email{loday@math.unistra.fr}
\maketitle

\date{\today}

\begin{abstract} We introduce a new one-person game similar to the Sudoku game. It is based on combinatorial objects called planar binary rooted trees. It is related to the four color conjecture. Its mathematical analysis makes use of the Tamari poset, hence the Stasheff associahedron.
\end{abstract}

\section{Introduction} 
Here is a first example of a YY game. In the following graph, label each edge either by $0$, $1$, or $2$, such that at each vertex, the three labels appear:

\medskip

{\centerline {\scalebox{0.6}{\includegraphics{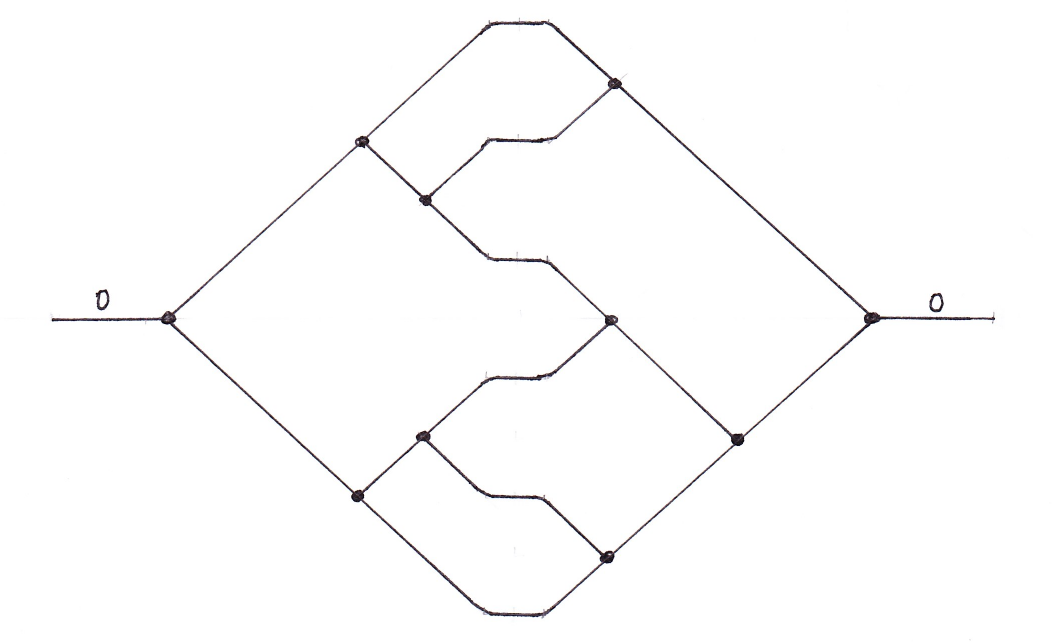}}}}

\medskip

A solution is to be found at the end of the paper.

\bigskip

We introduce a new one-person game based on planar binary trees. The aim is to label the edges of certain graphs, called YY-graphs (there are infinitely many of them), according to a simple rule. It is conjectured that for any one of these YY-graphs there is a solution. This conjecture is similar to a conjecture of Louis Kaufmann, whose proof would lead to a proof of the four colour theorem.

Some of the YY-graphs can be decomposed into smaller YY-graphs. This decomposition can be interpreted in terms of a certain poset structure on the set of planar binary trees, called the Tamari lattice. The geometric interpretation involves the Stasheff associahedron.

\section{Planar binary trees and YY-graphs} 
Let $PBT_{n}$ be the set of planar binary rooted trees with $n$ leaves:
\begin{displaymath}
PBT_1 = \{\  \vert \ \} , PBT_2=\big\{\arbreA \big\} , 
{PBT_3=\Big\{\arbreAB}\ ,\  {\arbreBA }\Big\}
\end{displaymath}
\begin{displaymath}
{PBT_4=\Bigg\{\arbreABC}\ , {\arbreBAC }, {\arbreACA }, {\arbreCAB },{\arbreCBA }\Bigg\}
\end{displaymath}

The integer $n$ is called the \emph{arity} of $t\in PBT_n$.


A planar binary rooted tree (or pb tree for short) has three different kinds of edges: the leaves (on top), the internal edges, and the root (on bottom). 

By definition a YY-graph is obtained by splicing two trees $s$ and $t$ with same arity along their respective leaves. We observe that a YY-graph $(s,t)$ is a graph whose vertices are ternary. 

\subsection{Example}\label{EX}

Let $s=\arbreABCplus$ and $t=\arbreCABplus$.

The YY-graph $(s,t)$ is
$$\graphUN$$
\subsection{The game} Starting with a YY-graph $(s,t)$ the game consists in labelling each edge by one of the three labels $0,1,2$ according to the following rules:

a) at each vertex we find the three labels, 

b) the labels of the two roots are the same (for instance $0$).

\medskip

A solution of the YY game presented in \ref{EX}  is 

$$\graphUNsol$$

The example of the introduction is made up out of the following trees $s$ and $t$ :

$s = \arbreBAEAB ,\qquad t = \arbreAEACA$

Check that a solution is given by the sequence $0, 1, 0, 1, 2, 1$ labelling the leaves (from left to right).

\section{Mathematical analysis of the YY-game}

\subsection{Labelling the leaves} The rule (a) implies that a labelling of a YY-graph is completely determined by the labelling of the edges which are the spliced leaves. Therefore the game consists in finding a labelling of the leaves such that, when we ``compute'' the labels for each tree we never meet the case where two inputs have the same label. The rule (b) says that the final value (the label of the root) is both time $0$. This  property can be handled mathematically under two different ways. One involves a magma, the other one a Lie algebra.

\subsection{Magma} We consider the following commutative magma. The underlying set is $X=\{0,1,2,\infty\}$. The binary operation is given by

\begin{displaymath}
\begin{array}{|c | c c c c |}
\hline
 & 0& 1&  2 & \infty \\
 \hline
0 & \infty & 2 & 1 & \infty \\
1 & 2 &\infty & 0 &  \infty \\
2 & 1 & 0 & \infty & \infty \\
\infty &\infty &\infty &\infty &\infty \\
\hline
\end{array}
\end{displaymath}

If $t$, resp.\ $s$, is a pb tree of arity $n$, any sequence $(x_1, \ldots, x_n)$ of elements in $X$ gives rise to an element in $X$, denoted
$t (x_1, \ldots, x_n)$, resp.\ $s (x_1, \ldots, x_n)$.  The YY-game is based on the following conjecture:

\medskip

\noindent {\bf Conjecture.} For any two planar binary trees $s$ and $t$ with same arity there is a sequence of elements $(x_1, \ldots, x_n)$ in $\{0,1,2\}$
such that 
$$t (x_1, \ldots, x_n)=s (x_1, \ldots, x_n)\neq \infty .$$

\subsection{Lie algebra}
Let us consider the set $Z=\{z_{0}, z_{1}, z_{2}\}$. We denote by $\mathbb{F}_{2}[X]$ the vector space spanned by $X$ over the field with $2$ elements. We equip it with the binary operation $[z_{i},z_{j}]$ given by the following table:

\begin{displaymath}
\begin{array}{|c | c c c |}
\hline
i\backslash j & z_{0} &  z_{1} &  z_{2}\\
\hline
z_{0} & 0 &  z_{2} &  z_{1}\\
z_{1} &  z_{2} & 0 &   z_{0}\\
z_{2} &  z_{1} &  z_{0}& 0 \\
\hline
\end{array}
\end{displaymath}

\begin{lemma}\label{lemmaLie} The operation $[-,-]$ is a Lie bracket on $\mathbb{F}_{2}[X]$.
\end{lemma}

\begin{proo} Since $+1=-1$ in $\mathbb{F}_{2}$ the operation is antisymmetric. The Jacobiator 
$$J(z_{i},z_{j},z_{k}):= [[z_{i},z_{j}], z_{k}]+[[z_{j},z_{k}], z_{i}]+[[z_{k},z_{i}], z_{j}]$$ 
takes the following value:

- if $i\neq j\neq k\neq i$, then we get
$$J(z_{0}, z_{1}, z_{2})= [z_{2}, z_{2}] + [z_{0}, z_{0}] + [ z_{1}, z_{1}]= 0,$$

- if $i=j\neq k$, then we get
$$ J(z_{0}, z_{0}, z_{1})= 0+[z_{2}, z_{0}] + [z_{2}, z_{0}]= 0,$$

- if $i=j=k$, then we get
$$J(z_{0}, z_{0}, z_{0})= 0+0+0=0.$$
\end{proo}

Any planar binary tree gives a way of computing an element in the Lie algebra out of the inputs. The conjecture can re-phrased as follows:

 For any two planar binary trees $s$ and $t$ with same arity there is a sequence of elements $(a_{1}, \ldots , a_{n})$ where $a_{i}$ is in $\{z_{0}, z_{1}, z_{2}\}$ such that
$$ s(a_{1}, \ldots , a_{n})= t(a_{1}, \ldots , a_{n})$$
and is not equal to $0$.

\subsection{The YY game, the Tamari lattice and the Stasheff polytope}

Let us suppose that we try to prove the conjecture by induction on the arity. Then it suffices to prove the conjecture for some of the pairs of trees in $PBT_{n}$, called \emph{prime YY-graphs} and defined as follows.

Given a YY-graph $(s,t)$, let us suppose that there is a sub-interval $[i, j]$ of the interval $[1,n]$ of the leaves, such that the sub-graph containing the leaves in this sub-interval form a sub-YY-graph. For instance, in the following example the interval $[2,4]$ has this property. 

\medskip

{\centerline {\scalebox{0.5}{\includegraphics{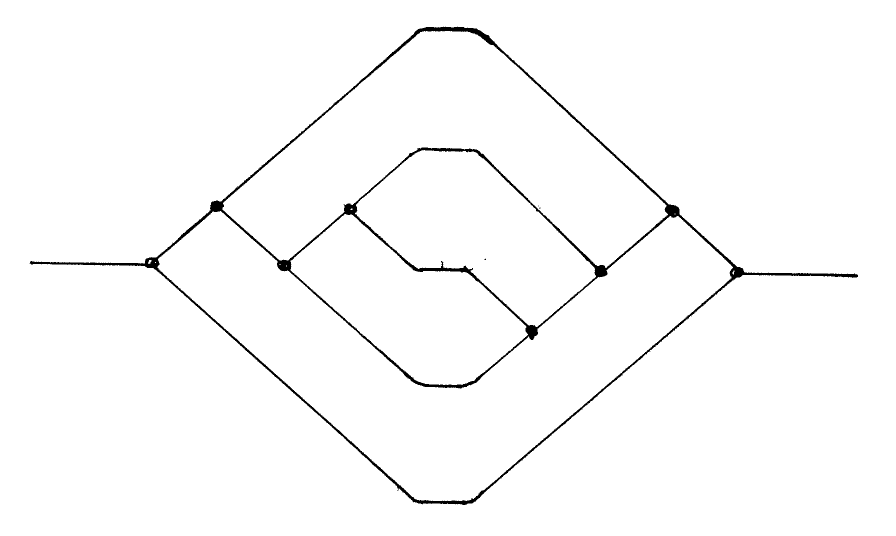}}}}

\medskip

If this sub-interval is strictly smaller than the full interval, then the game can be split into two games involving trees with strictly smaller arities. In our example they are the graphs:

\medskip

{\centerline {\scalebox{0.5}{\includegraphics{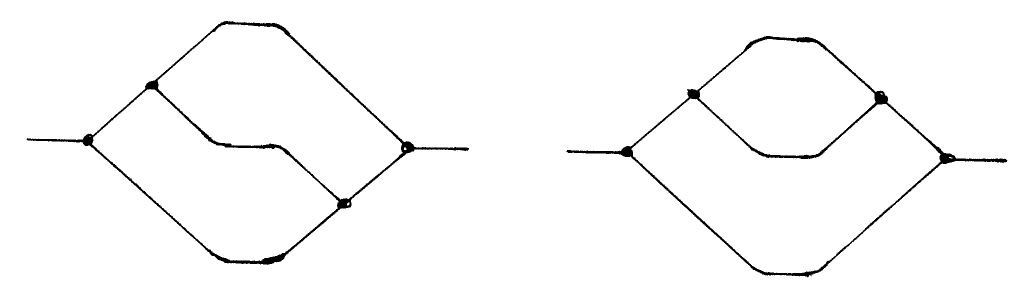}}}}

\medskip

On the left-hand side is the internal YY game, on the right-hand side is what is left when the internal YY game has been removed.

When such a sub-interval exists, we say that the graph is \emph{decomposable}. When there is no such sub-interval, then the YY-graph is said to be \emph{prime}. For instance the YY-graph of the introduction is prime.

There is a poset structure on the set $PBT_{n}$ which makes it into a lattice. It is called the Tamari order. The smallest, resp.\ largest, element for this order is the left comb denoted by $LC$, resp. right comb denoted by $RC$. Given two trees $s$ and $t$ in $PBT_{n}$, they admit a meet $t\wedge s$ and a join $t\vee s$. If we have 
$$t\wedge s= LC\quad \textrm{and} \quad t\vee s= RC$$
then the graph $(s,t)$ is prime, otherwise it is decomposable.

The pb trees can be interpreted as the vertices of a polytope called the Stasheff polytope, or the associahedron (cf.\ for instance \cite{JLLsta} for a simple explicit construction of it). More precisely each cell of this cellular complex can be encoded by a planar tree so that, if a cell $x$ belongs to the boundary of the cell $y$, then the tree encoding $y$ can be obtained from the tree $t$ encoding $x$ by shrinking some internal edges of $t$ to a point. The corolla (a tree with only one vertex)  is the tree encoding the big cell. The vertices of the Stasheff polytope are coded by the planar \emph{binary} trees.

If it happens that the two pb trees $s$ and $t$ (with same arity) encode vertices which belong to a subcell of the Stasheff polytope, then the YY-graph $(s,t)$ is decomposable.

\section{Relationship with the four color conjecture} In \cite{Kauffman} Louis Kauffman showed that the four color conjecture (4CC) follows from a conjecture which involves the following nonassociative algebra $\gg$ over a field $\KK$. As a vector space it is spanned by the elements $i,j,k$. The products are given by $[i,i]=[j,j]=[k,k]=0$ and $[i,j]=k=-[j,i], [j,k]=i=-[k,j], [k,i]=j=-[i,k]$. It is obviously nonassociative, however one can show that it is a Lie algebra. 

\begin{lemma} The product $[-,-]$ defined on $\gg$ is a Lie bracket. The Lie algebra $\gg$ can be presented with the generators $\{i,j\}$ and relations:
$$[[i,j],i]= j,\quad  [[j,i],j]= i .$$
\end{lemma}
\begin{proo}
Indeed, first the product is antisymmetric by definition. Second the Jacobi identity holds by a direct calculation as in the proof of Lemma \ref{lemmaLie}, for instance:
$$J(i,i,j) =[0,j]+ [k,i] + [-k,i] + = 0.$$

The verification of the second assertion is straightforward.
\end{proo}

Kauffman theorem asserts that, if, for any two ways $t$ and $s$ of parenthisizing a word $z_{1} \ldots z_{n}$, there is a nonzero solution to the equation $t(z_{1} \ldots z_{n})=s(z_{1} \ldots z_{n})$ in $\gg$, then the 4CC is true.

The difference between the aforementioned conjecture and Kaufman conjecture is, first, we forget about the signs, second we use only the generators and not the linear combinations.

\section{Appendix} A solution of the YY game 
mentioned in the introduction.

\medskip

{\centerline {\scalebox{0.6}{\includegraphics{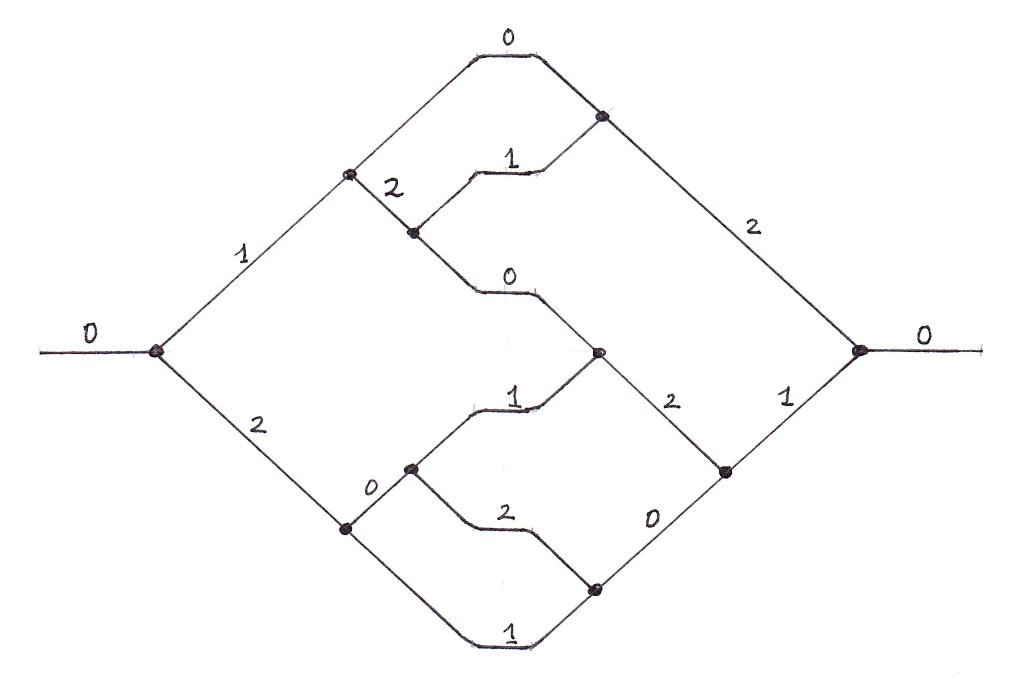}}}}

\end{document}